\pdfoutput=1
\documentclass{amsart}

\usepackage[T1]{fontenc}      
\usepackage{lmodern}          
\usepackage{microtype}        

\usepackage{amsmath, amssymb, amsthm}

\usepackage{graphicx}         
\usepackage{xcolor}           
\usepackage[margin=1.5in]{geometry} 
\usepackage{enumitem}         
\usepackage{setspace}
\usepackage{subcaption}

\usepackage[
    colorlinks=true,
    linkcolor=blue,
    citecolor=blue,
    urlcolor=blue
]{hyperref}

\theoremstyle{plain} 
\newtheorem{theorem}{Theorem}[section] 
\newtheorem{lemma}[theorem]{Lemma}         
\newtheorem{proposition}[theorem]{Proposition}

\theoremstyle{definition}
\newtheorem{definition}[theorem]{Definition}
\newtheorem{example}[theorem]{Example}

\theoremstyle{remark}
\newtheorem{remark}[theorem]{Remark}

\title[REGION COLORING OF SPATIAL GRAPHS]{Region Coloring of Unoriented Spatial Graphs using Knot-theoretic Ternary Quasigroups}
\author{Mitsuru Kuramoto}
\address{Department of Mathematics, Graduate School of Science, The University of Osaka Toyonaka, Osaka 560-0043, JAPAN}
\email{u119849e@ecs.osaka-u.ac.jp}

\begin{document}
\onehalfspacing 

\begin{abstract}
    We construct coloring conditions around vertices for region colorings of unoriented spatial graph diagrams by knot-theoretic ternary quasigroups.
    These conditions require cyclic products of left and right operators around each vertex to be the identity.
    We prove that the resulting coloring sets are preserved under the generalized Reidemeister moves and relate our construction to the palette framework.
    We show that the vertex coloring conditions are automatic for Dehn tribrackets and determine the associated palettes for unoriented Alexander tribrackets.
    We also introduce an additional restriction on vertex colorings that is preserved under the fifth Reidemeister move for every unoriented knot-theoretic ternary quasigroup.
    This restriction gives rise to further coloring invariants for Dehn tribrackets and unoriented Alexander tribrackets.
\end{abstract}

\maketitle

\section{introduction}
A knot-theoretic ternary quasigroup (KTQ) is an algebraic structure equipped with a ternary operation satisfying certain axioms.
The use of ternary operations in knot theory was introduced by Niebrzydowski in \cite{MR3205573},
and was subsequently developed into the theory of KTQs in 
\cite{niebrzydowski2018ternaryquasigroupsknottheory}.
This ternary operation is often referred to as a Niebrzydowski tribracket, or simply a tribracket.
This algebra is derived from the structure of Reidemeister moves on link diagrams.
It provides a framework for coloring the regions of a diagram,
where the number of such colorings yields a link invariant.

The theory was further extended through homological constructions
and cocycle invariants for knots and knotted surfaces in \cite{MR4171568}.

A spatial graph is a graph embedded in $\mathbb{R}^{3}$, 
serving as a natural generalization of knots and links. 
When extending KTQ colorings to spatial graphs, 
the presence of vertices poses a significant topological hurdle. 
In \cite{MR3900876}, a KTQ coloring for Y-oriented trivalent spatial graphs was established by defining a local algebraic structure around the vertices. 
More recently, Oshiro and Oyamaguchi \cite{MR4259504} introduced the concept of \textit{palettes} of KTQ colorings for spatial graphs, which are local coloring conditions around a vertex preserved under the generalized Reidemeister moves, 
and they completely classified the palettes for Dehn $p$-colorings.

In this paper, we construct the largest unoriented palette
\[
\mathcal{P}(X)=\bigcup_{n\ge 3}\mathcal{P}_{(n)}
\]
associated with an arbitrary unoriented KTQ.
The set $\mathcal{P}_{(n)}$ consists of admissible color tuples at an
$n$-valent vertex.
To achieve this, we introduce right and left operators, $R$ and $L$,
and formulate general coloring conditions around an $n$-valent vertex.

We show that, for Dehn tribrackets,
every tuple satisfies the vertex coloring condition.
For unoriented Alexander tribrackets, we determine
the associated palettes.

Furthermore, we introduce an additional subset
\[
\mathcal{Q}_{(n)}(X,[\cdot,\cdot,\cdot]) \subset X^n,
\]
which provides a further restriction on vertex colorings,
and prove that the resulting restricted colorings are preserved
under all generalized Reidemeister moves for Dehn tribrackets and
unoriented Alexander tribrackets.

The remainder of this paper is organized as follows.
In Section 2, we recall the definition of an unoriented KTQ and the
region coloring conditions for classical links.
In Section 3, we define two bijective operators $R$ and $L$ on a KTQ
and verify several algebraic identities involving these operators.
In Section 4, we construct the vertex conditions and the associated
largest palette.
We also give a finite example and show that the vertex conditions are
automatic for Dehn tribrackets.
For unoriented Alexander tribrackets, we completely determine the
associated palettes.
Finally, in Section 5, we introduce an additional coloring restriction
at vertices and the associated subsets $\mathcal{Q}_{(n)}$.
We prove that the resulting restricted colorings are preserved under
RV moves for every unoriented KTQ, and under all generalized
Reidemeister moves for Dehn tribrackets and for unoriented Alexander
tribrackets.

\section{knot-theoretic ternary quasigroups}
In this section, we review the concept of KTQ originally introduced for oriented links, and then discuss region colorings of unoriented link diagrams.
We also present several specific examples and describe their properties
\cite{MR4547475, MR3900876, MR4259504}.

In this paper, we define a modified KTQ, which is a pair of a set $X$ and a ternary operation $[\cdot,\cdot,\cdot]$.

\begin{definition}\label{KTQdef}
    Let $X$ be a set. A ternary operation 
    $[\cdot, \cdot , \cdot]:X^3 \to X$ is called a \textit{tribracket} if it satisfies the following conditions:
    
    (KTQ1) For all $x, y, z \in X$,
        \begin{enumerate}
            \item there exists a unique element $w \in X$ such that $[w, x, y] = z$,
            \item there exists a unique element $w \in X$ such that $[x, w, y] = z$,
            \item there exists a unique element $w \in X$ such that $[x, y, w] = z$.
        \end{enumerate}  
        
    (KTQ2) For all $x, y, z, w \in X$, the following equations hold:
    \[
    [[y,x,z],y,[y,x,w]] = [[y,x,z],z,[z,x,w]] = [[y,x,w],w,[z,x,w]].
    \]
    We call the pair $(X, [\cdot,\cdot,\cdot])$ a \textit{Knot-theoretic Ternary Quasigroup} or \textit{KTQ}.
\end{definition}

\begin{remark}\label{Remark-KTQdef}
    In the traditional convention (e.g., \cite{MR4259504, MR3900876}), a KTQ is defined as a pair $(X, [\cdot,\cdot,\cdot]_{0})$ consisting of a set $X$ and a ternary operation $[\cdot, \cdot , \cdot]_{0}$ that satisfies the same invertibility conditions as (KTQ1) and the following equation instead of (KTQ2):
    \[
    [y,[x,y,z]_{0},[x,y,w]_{0}]_{0} = [z,[x,y,z]_{0},[x,z,w]_{0}]_{0} = [w,[x,y,w]_{0},[x,z,w]_{0}]_{0}
    \]
    for all $x, y, z, w \in X$. 
    Our modified operation $[\cdot,\cdot,\cdot]$ is related to the traditional one via $[x,y,z] = [y,x,z]_{0}$, 
    which corresponds to a permutation of the first two arguments.
\end{remark}

\begin{example}\label{KTQdef-dehn}
    Let $X$ be a group, and define a ternary operation by
    \[
    [x,y,z] = x y^{-1}z.
    \] 
    This ternary operation is called the \textit{Dehn tribracket}.
\end{example}

\begin{example}\label{KTQdef-alex}
    Let $X$ be a commutative ring with identity, and let $t$ and $s$ be units in $X$. 
    Define a ternary operation on $X$ by
    \[
    [x,y,z] = tx - tsy + sz. 
    \]
    A KTQ equipped with this operation is called an \textit{Alexander tribracket}.
    In particular, if $X = \mathbb{Z}/p\mathbb{Z}$ for an integer $p \geq 2$, 
    and we set $t = -1$ and $s = 1$, 
    we obtain the special \textit{Alexander tribracket}. 
    This structure is related to \textit{Dehn $p$-coloring}, 
    which is a coloring system for the regions of a knot diagram.
\end{example}


KTQs were introduced to construct oriented knot invariants 
via colorings of the regions of a knot diagram. 
In this part, we review the notion of an \textit{unoriented tribracket}, 
which is a special property concerning colorings of unoriented diagrams.

\begin{definition}\label{unoridef}
    Let $(X, [\cdot, \cdot, \cdot])$ be a KTQ. 
    This tribracket is called \textit{unoriented} if it satisfies the following conditions: \\
    (KTQ3) For all $x, y, z \in X$,
    \[
        (1)\quad z=[x,[y,z,x],y],\quad
        (2)\quad z=[[z,y,x],x,y], \quad
        (3)\quad z=[y,x,[x,y,z]].
    \]
    If a tribracket $[\cdot,\cdot,\cdot]$ of a KTQ $(X,[\cdot,\cdot,\cdot])$ is unoriented, we also call $(X,[\cdot,\cdot,\cdot])$ an \textit{unoriented KTQ}. 
\end{definition}

\begin{remark}\label{remark-unoridef}
    Any two of the conditions (1), (2) and (3) in (KTQ3) imply the remaining one. 
    For instance, assuming (2) and (3), we can derive (1) as follows :
    setting $w=[y,z,x]$, condition (3) gives $[z,y,w]=x$, and condition (2) gives $[x,[y,z,x],y]=[[z,y,w],w,y]=z$. 
\end{remark}
\begin{remark}\label{previous-unori-def}
    Let $(X,[\cdot,\cdot,\cdot])$ be a KTQ in our convention and $(X,[\cdot,\cdot,\cdot]_{0})$ be a KTQ in the sense of previous studies such that $[a,b,c]=[b,a,c]_{0}$.
    In Definition 8.1 of \cite{MR4259504}, an unoriented tribracket $[\cdot,\cdot,\cdot]_{0}$ is defined by the conditions:\\
    $(\text{KTQ3})_{0}$ For all $x, y, z \in X$,
        \[
        z=[x,y,[y,x,z]_{0}]_{0}\quad \text{and} \quad z=[x,[y,z,x]_{0},y]_{0}.
        \]
    Rewriting these equations via $[x,y,z]=[y,x,z]_{0}$ and applying Remark~\ref{remark-unoridef},
    we see that our definition of an unoriented tribracket in Definition~\ref{unoridef} is equivalent to that in \cite{MR4259504}.
\end{remark}

By direct calculation, 
we can see that several tribrackets are \textit{unoriented}.

\begin{example}\label{unori-dehn}
    The Dehn tribracket $[x,y,z] = xy^{-1}z$ is \textit{unoriented}.
\end{example}

\begin{example}\label{unori-alex}
    The Alexander tribracket $[x,y,z] = tx - tsy + sz$ is \textit{unoriented} if the units $t$ and $s$ satisfy $t^2 = s^2 = 1$.
\end{example}

\begin{figure}[htbp]
    \centering
    \begin{minipage}[c]{0.45\textwidth}
        \centering
        \includegraphics{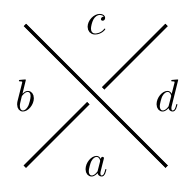}
    \end{minipage}%
    \begin{minipage}[c]{0.45\textwidth}
        \centering
        $\begin{aligned}
            \mathrm{(i)}   & \quad a=[d,c,b], \\
            \mathrm{(ii)}  & \quad b=[c,d,a], \\
            \mathrm{(iii)} & \quad c=[b,a,d], \\
            \mathrm{(iv)}  & \quad d=[a,b,c].
        \end{aligned}$
    \end{minipage}%
    \caption{Coloring condition for an unoriented diagram}
    \label{fig:coloring_condition}
\end{figure}

We now describe the coloring conditions for unoriented knot diagrams using an \textit{unoriented tribracket}. 
Let $(X, [\cdot, \cdot, \cdot])$ be a KTQ equipped with an \textit{unoriented tribracket}, 
let $D$ be a knot or link diagram, 
and let $R(D)$ be the set of regions of $D$. 
A map $c: R(D) \to X$ is called a \textit{coloring map} if it satisfies the coloring condition depicted in Figure~\ref{fig:coloring_condition} at each crossing.

To determine the color of a region at a crossing, 
read the colors of the other three regions in the following order: 
across the over-arc, diagonally opposite, and across the under-arc. 
The color of the chosen region is the tribracket of these three colors, in this order. 
The four equations (i)-(iv) are interpreted as follows:
\[
    \mathrm{(i)'} \quad a = [c,d,b]_{0}, \quad \mathrm{(ii)'} \quad b=[d,c,a]_{0}, \quad \mathrm{(iii)'} \quad c=[a,b,d]_{0}, \quad \text{and} \quad \mathrm{(iv)'} \quad d=[b,a,c]_{0}.
\]
Furthermore, the four equations are mutually equivalent \cite{MR4259504}. 
Thus, satisfying any single condition is sufficient for the coloring condition.

\begin{figure}[htbp]
    \centering
    
    \begin{minipage}{\textwidth}
        \centering
        \includegraphics{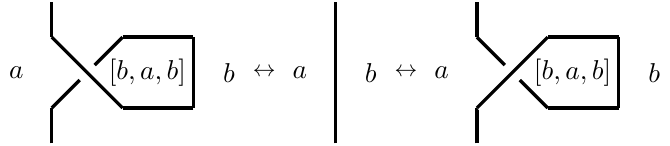}
        \\[0.5em]
        {\small (a) Reidemeister move I (RI)}
    \end{minipage}

    \vspace{1.5em}

    \begin{minipage}{\textwidth}
        \centering
        \includegraphics{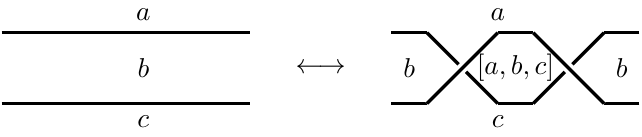}
        \\[0.5em]
        {\small (b) Reidemeister move II (RII)}
    \end{minipage}

    \vspace{1.5em}

    \begin{minipage}{\textwidth}
        \centering
        \includegraphics[width=\linewidth]{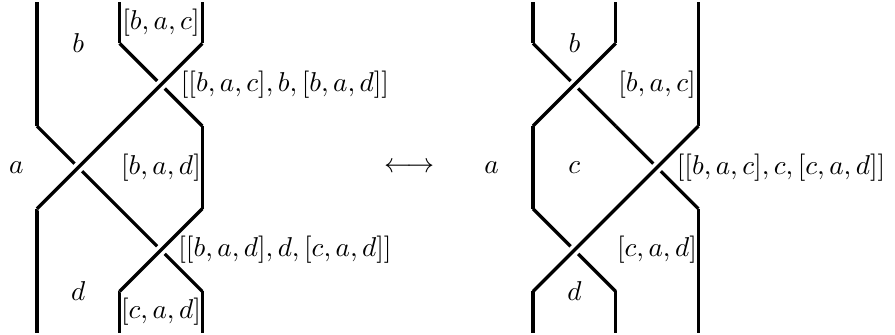}
        \\[0.5em]
        {\small (c) Reidemeister move III (RIII)}
    \end{minipage}

    \vspace{1em}
    \caption{Coloring preservation under Reidemeister moves}
    \label{fig:colored_reidemeister_moves}
\end{figure}

We denote by $\mathrm{Col}_X(D)$ the set of coloring maps. 
As in the case of quandle colorings, there exists a bijection between the sets of colorings under each Reidemeister move depicted in Figure~\ref{fig:colored_reidemeister_moves}.


\begin{proposition}\label{prop-invariance-classical}
    Let $(X,[\cdot,\cdot,\cdot])$ be an unoriented KTQ. 
    Let $D$ and $D'$ be diagrams of a link $L$. 
    Then, there exists a bijection between $\mathrm{Col}_X(D)$ and $\mathrm{Col}_X(D')$.
    Consequently, if $|X|<\infty$, the number of colorings $|\mathrm{Col}_X(D)|$ is an invariant of unoriented links.
\end{proposition}


\section{operator Identities for unoriented KTQs}
Although a tribracket is a ternary operation, fixing two of its elements yields a bijection from $X$ to $X$. 
Each of these bijections acts as an operator. 
In this section, we introduce the operator identities for unoriented KTQs.

\begin{definition}\label{def-RLope}
    Let $(X, [\cdot, \cdot, \cdot])$ be an unoriented KTQ.
    For any $a, b \in X$, we define two maps from $X$ to $X$ by
    \begin{enumerate}[label=\upshape(\arabic*)]
        \item $R(a,b) \colon X \to X, \quad x \mapsto [x,a,b] =: (x)R(a,b)$,
        \item $L(b,a) \colon X \to X, \quad x \mapsto [b,a,x] =: L(b,a)(x)$.
    \end{enumerate}
\end{definition}

\begin{definition}\label{def-opcomp}
    Let $(X, [\cdot ,\cdot ,\cdot])$ be an unoriented KTQ.
    For any $a, b, c, d \in X$, the compositions of these operators are denoted by
    \begin{enumerate}[label=\upshape(\arabic*)]
        \item $(x)R(a,b)R(c,d) := [[x,a,b],c,d]$,
        \item $L(d,c)L(b,a)(x) := [d,c,[b,a,x]]$.
    \end{enumerate}
\end{definition}

From the condition (KTQ3) of an unoriented tribracket, we obtain the following relations:
\begin{enumerate}[label=\upshape(\arabic*)]
    \item $R(a,b)R(b,a)=\mathrm{id}_{X}$,
    \item $L(a,b)L(b,a)=\mathrm{id}_{X}$.
\end{enumerate}


Let $G_{R}$ be the group generated by $\left\{ R(x,y) \mid x,y\in X \right\}$ with the product defined as in Definition~\ref{def-opcomp}~(1), 
and let $G_{L}$ be the group generated by $\left\{ L(y,x) \mid x,y\in X \right\}$ with the product defined as in Definition~\ref{def-opcomp}~(2). 
Thus, $G_{R}$ acts on $X$ from the right, and $G_{L}$ acts on $X$ from the left. 
We now present several identities among these operators.

\begin{proposition}\label{prop-RL-identities}
    Let $(X, [\cdot ,\cdot ,\cdot])$ be an unoriented KTQ.
    For all $a, b, c, x, y \in X$, the following relations among the operators hold:
    \begin{enumerate}[label=\upshape(\arabic*)]
        \item $R(a,b)R(b,c) = R(a,[c,b,a])R([c,b,a],c)$,
        \item $R(a,b)R(b,c) = R(a,[a,b,c])R([a,b,c],c)$,
        \item $L(c,b)L(b,a) = L(c,[c,b,a])L([c,b,a],a)$,
        \item $L(c,b)L(b,a) = L(c,[a,b,c])L([a,b,c],a)$,
        \item If $y=L(b,a)(x)$, then $L(y,x) = L(y,b)L(b,a)L(a,x)$,
        \item If $y=L(b,a)(x)$, then $R(x,y) = R(x,a)R(a,b)R(b,y)$,
        \item If $y=(x)R(a,b)$, then $L(y,x) = L(y,b)L(b,a)L(a,x)$,
        \item If $y=(x)R(a,b)$, then $R(x,y) = R(x,a)R(a,b)R(b,y)$.
    \end{enumerate}
\end{proposition}


\begin{proof}
    These relations follow from the axioms (KTQ2) and (KTQ3) by direct computation.
    Recall that for all $x, y, z, w \in X$, (KTQ2) states
    \[
        [[y,x,z],y,[y,x,w]] = [[y,x,z],z,[z,x,w]] = [[y,x,w],w,[z,x,w]],
    \]
    and (KTQ3) states that for all $x, y, z \in X$,
    \[
        (1)\quad z=[x,[y,z,x],y],\quad
        (2)\quad z=[[z,y,x],x,y], \quad
        (3)\quad z=[y,x,[x,y,z]].
    \]

    \noindent (1) For any $x \in X$,
    \begin{align*}
        (x)R(a,b)R(b,c)
        &= [[x,a,b],b,c] \\
        &= [[x,a,b],b,[[c,b,a],a,b]] \\
        &= [[x,a,[c,b,a]],[c,b,a],[[c,b,a],a,b]] \\
        &= [[x,a,[c,b,a]],[c,b,a],c] \\
        &= (x)R(a,[c,b,a])R([c,b,a],c).
    \end{align*}

    \noindent (2) For any $x \in X$,
    \begin{align*}
        (x)R(a,b)R(b,c)
        &= [[x,a,b],b,c] \\
        &= [[x,a,b],b,[b,a,[a,b,c]]] \\
        &= [[x,a,[a,b,c]],[a,b,c],[b,a,[a,b,c]]] \\
        &= [[x,a,[a,b,c]],[a,b,c],c] \\
        &= (x)R(a,[a,b,c])R([a,b,c],c).
    \end{align*}

    \noindent (3) For any $x \in X$,
    \begin{align*}
        L(c,b)L(b,a)(x)
        &= [c,b,[b,a,x]] \\
        &= [[[c,b,a],a,b],b,[b,a,x]] \\
        &= [[[c,b,a],a,b],[c,b,a],[[c,b,a],a,x]] \\
        &= [c,[c,b,a],[[c,b,a],a,x]] \\
        &= L(c,[c,b,a])L([c,b,a],a)(x).
    \end{align*}

    \noindent (4) For any $x \in X$,
    \begin{align*}
        L(c,b)L(b,a)(x)
        &= [c,b,[b,a,x]] \\
        &= [[b,a,[a,b,c]],b,[b,a,x]] \\
        &= [[b,a,[a,b,c]],[a,b,c],[[a,b,c],a,x]] \\
        &= [c,[a,b,c],[[a,b,c],a,x]] \\
        &= L(c,[a,b,c])L([a,b,c],a)(x).
    \end{align*}

    \noindent (5) Suppose $y=L(b,a)(x)$. For any $z \in X$,
    \begin{align*}
        L(y,x)(z)
        &= L(L(b,a)(x), x)(z) \\
        &= [[b,a,x],x,z] \\
        &= [[b,a,x],x,[x,a,[a,x,z]]] \\
        &= [[b,a,x],b,[b,a,[a,x,z]]] \\
        &= [y,b,L(b,a)L(a,x)(z)] \\
        &= L(y,b)L(b,a)L(a,x)(z).
    \end{align*}

    \noindent (6) Suppose $y=L(b,a)(x)$. For any $z \in X$,
    \begin{align*}
        (z)R(x,y)
        &= (z)R(x,L(b,a)(x)) \\
        &= [z,x,[b,a,x]] \\
        &= [[[z,x,a],a,x],x,[b,a,x]] \\
        &= [[[z,x,a],a,b],b,[b,a,x]] \\
        &= [[[z,x,a],a,b],b,y] \\
        &= (z)R(x,a)R(a,b)R(b,y).
    \end{align*}

    \noindent (7) Suppose $y=(x)R(a,b)$. For any $z \in X$,
    \begin{align*}
        L(y,x)(z)
        &= L((x)R(a,b), x)(z) \\
        &= [[x,a,b],x,z] \\
        &= [[x,a,b],x,[x,a,[a,x,z]]] \\
        &= [[x,a,b],b,[b,a,[a,x,z]]] \\
        &= [y,b,L(b,a)L(a,x)(z)] \\
        &= L(y,b)L(b,a)L(a,x)(z).
    \end{align*}

    \noindent (8) Suppose $y=(x)R(a,b)$. For any $z \in X$,
    \begin{align*}
        (z)R(x,y)
        &= (z)R(x,(x)R(a,b)) \\
        &= [z,x,[x,a,b]] \\
        &= [[[z,x,a],a,x],x,[x,a,b]] \\
        &= [[[z,x,a],a,b],b,[x,a,b]] \\
        &= [[[z,x,a],a,b],b,y] \\
        &= (z)R(x,a)R(a,b)R(b,y).
    \end{align*}
\end{proof}


\section{coloring conditions for spatial graphs}
In this paper, we assume that a finite graph may contain loops, multiple edges, and circle components without vertices, 
and that the degree of every vertex is at least three.
A tame embedding of such a finite graph into $\mathbb{R}^3$ is called a \textit{spatial graph}. 
Two spatial graphs are \textit{equivalent} (or \textit{ambient isotopic}) if one can be transformed into the other by an ambient isotopy of $\mathbb{R}^3$.
A diagram of a spatial graph is a regular projection to $\mathbb{R}^2$ equipped with over-under information at each crossing, 
where the multiple points of the edges are finitely many transverse double points.
The \textit{regions} of a spatial graph diagram $D$ are the connected components in the complement of its projection in $\mathbb{R}^2$,
and we denote this set of regions by $R(D)$.

\begin{figure}[htbp]
    \centering
    
    \begin{minipage}{0.48\textwidth}
        \centering
        \includegraphics{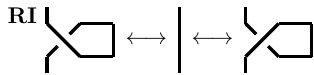}
    \end{minipage}%
    \hfill
    \begin{minipage}{0.48\textwidth}
        \centering
        \includegraphics{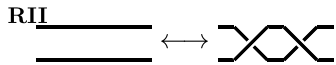}
    \end{minipage}
    
    \vspace{1.5em} 
    
    \begin{minipage}{0.48\textwidth}
        \centering
        \includegraphics{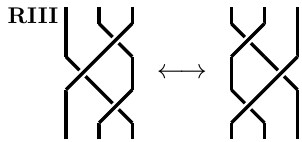}
    \end{minipage}%
    \hfill
    \begin{minipage}{0.48\textwidth}
        \centering
        \includegraphics{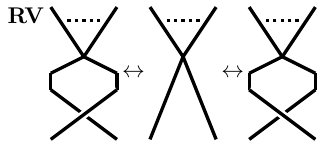}
    \end{minipage}
    
    \vspace{1.5em} 

    \begin{minipage}{\textwidth}
        \centering
        \includegraphics{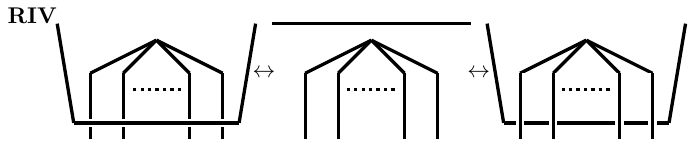}
    \end{minipage}

    \caption{Generalized Reidemeister moves for spatial graphs (RI--RV).}
    \label{fig:reidemeister_moves}
\end{figure}

\begin{lemma}[\cite{MR946218}]\label{thm-reidemeister}
    Let $D_1$ and $D_2$ be two diagrams representing equivalent spatial graphs. 
    Then $D_1$ and $D_2$ are related by a finite sequence of planar isotopies and generalized Reidemeister moves (RI--RV) as shown in Figure~\ref{fig:reidemeister_moves}.
\end{lemma}


When extending KTQ colorings of knots to spatial graphs, 
the primary challenge lies in formulating an appropriate coloring condition around vertices. 
In this section, we introduce a KTQ coloring condition for spatial graph diagrams, 
and show that the resulting set of colorings yields an invariant of spatial graphs. 
Furthermore, we show that every tuple satisfies the vertex condition for a Dehn tribracket. 
For an unoriented Alexander tribracket, we determine the associated palettes.

We define a KTQ coloring condition for spatial graphs.

\begin{definition}\label{def-RL-id}
    Let $(X,[\cdot,\cdot,\cdot])$ be an unoriented KTQ, and let $n\geq 3$ be an integer.
    For an $n$-tuple $(x_{1}, \dots, x_{n})\in X^n$, we consider the following two conditions:
    \begin{align}
        R(x_1,x_2)R(x_2,x_3) \cdots R(x_{n-1},x_n)R(x_n,x_1) &= \mathrm{id}_X,\label{eq:R-id}\\
        L(x_1,x_n)L(x_n,x_{n-1}) \cdots L(x_3,x_2)L(x_2,x_1) &= \mathrm{id}_X.\label{eq:L-id}
    \end{align}
    We denote by $\mathcal{R}_{(n)}(X,[\cdot,\cdot,\cdot])$ and $\mathcal{L}_{(n)}(X,[\cdot,\cdot,\cdot])$ the subsets of $X^n$ consisting of $n$-tuples satisfying \eqref{eq:R-id} and \eqref{eq:L-id}, respectively.
    We define
    \[
    \mathcal{P}_{(n)}(X,[\cdot,\cdot,\cdot]) := \mathcal{R}_{(n)}(X,[\cdot,\cdot,\cdot])\cap \mathcal{L}_{(n)}(X,[\cdot,\cdot,\cdot]).
    \]
    When $(X,[\cdot,\cdot,\cdot])$ is clear from the context, we simply write $\mathcal{R}_{(n)}$, $\mathcal{L}_{(n)}$, and $\mathcal{P}_{(n)}$.
\end{definition}

\begin{figure}[htbp]
    \centering
    \begin{minipage}[c]{0.45\textwidth}
        \centering
        \includegraphics{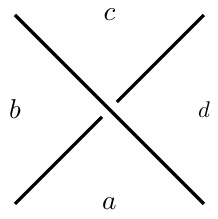}
    \end{minipage}\hfill
    \begin{minipage}[c]{0.45\textwidth}
        \centering
        \includegraphics{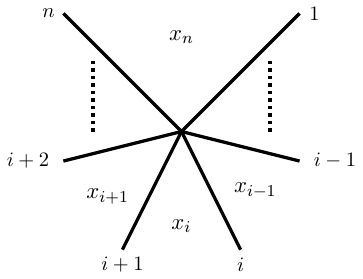}
    \end{minipage}
    \vspace{0.8em}
    \caption{Coloring condition at a crossing and an $n$-valent vertex}
    \label{fig:vertex_coloring_condition}
\end{figure}

\begin{definition}\label{def-col-condi}
    Let $D$ be a spatial graph diagram in $\mathbb{R}^2$ and $(X, [\cdot , \cdot, \cdot])$ be an unoriented KTQ. 
    Let $R(D)$ be the set of regions of $D$.
    A map $c : R(D) \to X$ is called a \textit{KTQ coloring} (or simply a \textit{coloring map}) of $D$ if it satisfies the following conditions:
    \begin{enumerate}
        \item \textbf{Coloring condition at crossings}\\
        At each crossing, the colors $a,b,c,d \in X$ assigned to the four surrounding regions satisfy $d=[a,b,c]$, as depicted on the left of Figure~\ref{fig:vertex_coloring_condition}.
        \item \textbf{Coloring condition at vertices}\\
        At each $n$-valent vertex, as depicted on the right of Figure~\ref{fig:vertex_coloring_condition}, read the colors of the $n$ local sectors in clockwise order to obtain an $n$-tuple $(x_{1},\dots , x_{n})$.
        We require that $(x_{1},\dots,x_{n})\in\mathcal{P}_{(n)}$.        
    \end{enumerate}
\end{definition}


\begin{remark}\label{remark-cyclic-ordering}
    Note that this coloring condition around a vertex is well-defined, being independent of the choice of the initial region in the cyclic ordering. 
    Since the operators belong to the group $G_R$ or $G_L$, 
    if a product of operators equals the identity $\mathrm{id}_X$, 
    any cyclic permutation of the product also equals $\mathrm{id}_X$ by conjugation.
    The condition is also independent of the direction of the cyclic ordering, 
    since reversing the order replaces each cyclic product by its inverse, up to a cyclic permutation.
\end{remark}

We denote by $\operatorname{Col}_X(D)$ the set of all coloring maps of an unoriented spatial graph diagram $D$ with values in $X$.

\begin{theorem}\label{mainresult-1}
    Let $D_1$ and $D_2$ be diagrams representing the same spatial graph, 
    and let $(X, [\cdot, \cdot, \cdot])$ be an unoriented KTQ.
    Then, there exists a bijection between $\operatorname{Col}_X(D_1)$ and $\operatorname{Col}_X(D_2)$.
\end{theorem}

\begin{proof}
    By Lemma~\ref{thm-reidemeister}, 
    $D_1$ and $D_2$ are related by a finite sequence of generalized Reidemeister moves (RI--RV). 
    By Proposition~\ref{prop-invariance-classical}, there exist one-to-one correspondences between the coloring sets under the standard Reidemeister moves RI, RII, and RIII. 
    Therefore, it suffices to verify the bijections for the spatial graph moves RIV and RV. 

    \noindent\textbf{RIV moves}

    We consider the case where a strand passes over $n$ strands as depicted in Figure~\ref{fig:strand_over}.

    \begin{figure}[htbp]
    \centering
    \includegraphics{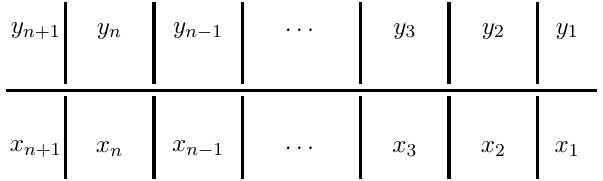}
    \caption{A strand passing over $n$ strands.}
    \label{fig:strand_over}
    \end{figure}

    The conditions at each crossing yield the following relations for $x_{1},\dots, x_{n+1}, y_{1}, \dots, y_{n+1}$:

    \noindent\textbf{Strand over-pass relations}
    \begin{align}
        y_{2} &= [x_{2},x_{1},y_{1}] = L(x_{2},x_{1})(y_{1}), \label{eq:over-first}\\
        y_{3} &= [x_{3},x_{2},y_{2}] = L(x_{3},x_{2})(y_{2})=L(x_{3},x_{2})L(x_{2},x_{1})(y_{1}), \\
        &\hspace{0.5em}\vdots\nonumber\\
        y_{n} &= [x_{n}, x_{n-1}, y_{n-1}] = L(x_{n},x_{n-1})(y_{n-1})= L(x_{n},x_{n-1})\cdots L(x_{2},x_{1})(y_{1}), \label{eq:over-last}\\
        y_{n+1} &= [x_{n+1}, x_{n}, y_{n}] = L(x_{n+1},x_{n})(y_{n})= L(x_{n+1},x_{n})\cdots L(x_{2},x_{1})(y_{1}), \label{eq:over-last2}\\
        x_{2} &= [y_{2},y_{1},x_{1}] = L(y_{2},y_{1})(x_{1}), \label{eq:over-x2}\\
        x_{3} &= [y_{3},y_{2},x_{2}] = L(y_{3},y_{2})(x_{2})= L(y_{3},y_{2})L(y_{2},y_{1})(x_{1}), \\
        &\hspace{0.5em}\vdots\nonumber\\
        x_{n} &= [y_{n},y_{n-1},x_{n-1}] = L(y_{n},y_{n-1})(x_{n-1})= L(y_{n},y_{n-1})\cdots L(y_{2},y_{1})(x_{1}), \\
        x_{n+1} &= [y_{n+1},y_{n},x_{n}] = L(y_{n+1},y_{n})(x_{n})= L(y_{n+1},y_{n})\cdots L(y_{2},y_{1})(x_{1})\label{eq:over-xn1}.
    \end{align}

    We consider the case where a strand passes under $n$ strands as depicted in Figure~\ref{fig:strand_under}.

    \begin{figure}[htbp]
    \centering
    \includegraphics{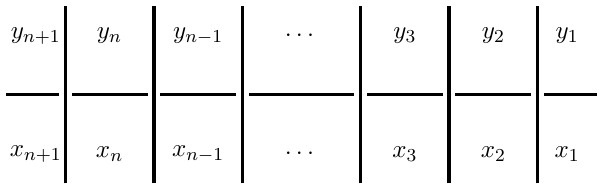}
    \caption{A strand passing under $n$ strands.}
    \label{fig:strand_under}
    \end{figure}

    Similarly, the coloring conditions at each crossing yield the following relations:


    \noindent\textbf{Strand under-pass relations}
    \begin{align}
        y_{2} &= [y_{1},x_{1},x_{2}] = (y_1)R(x_{1},x_{2}), \label{eq:under-first}\\
        y_{3} &= [y_{2},x_{2},x_{3}] = (y_{2})R(x_{2},x_{3}) = (y_1)R(x_{1},x_{2})R(x_{2},x_{3}), \\
        &\hspace{0.5em}\vdots\nonumber\\
        y_{n} &= [y_{n-1},x_{n-1},x_{n}] = (y_{n-1})R(x_{n-1},x_{n}) = (y_1)R(x_{1},x_{2})\cdots R(x_{n-1},x_{n}), \\
        y_{n+1} &= [y_{n},x_{n},x_{n+1}] = (y_{n})R(x_{n},x_{n+1}) = (y_1)R(x_{1},x_{2})\cdots R(x_{n},x_{n+1}), \\
        x_{2} &= [x_{1},y_{1},y_{2}] = (x_{1})R(y_{1},y_{2}), \\
        x_{3} &= [x_{2},y_{2},y_{3}] = (x_{2})R(y_{2},y_{3}) = (x_{1})R(y_{1},y_{2})R(y_{2},y_{3}), \\
        &\hspace{0.5em}\vdots\nonumber\\
        x_{n} &= [x_{n-1},y_{n-1},y_{n}] = (x_{n-1})R(y_{n-1},y_{n}) = (x_{1})R(y_{1},y_{2})\cdots R(y_{n-1},y_{n}), \\
        x_{n+1} &= [x_{n},y_{n},y_{n+1}] = (x_{n})R(y_{n},y_{n+1}) = (x_{1})R(y_{1},y_{2})\cdots R(y_{n},y_{n+1}). \label{eq:under-last}
    \end{align}

\begin{figure}[htbp]
        \centering
        \includegraphics{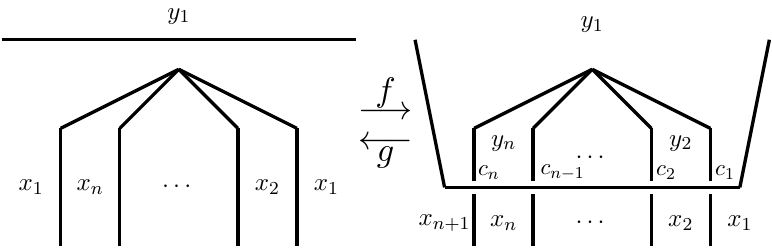}
        \caption{Colorings under the spatial graph move RIV (the strand passes over the vertex).}
        \label{fig:strand_over_vertex}
    \end{figure}
    
    Considering the local coloring correspondence $f$ depicted in Figure~\ref{fig:strand_over_vertex}, such that $x_{n+1}=x_{1}$, 
    the values of the new regions $y_{2}, y_{3}, \dots, y_{n}$ are uniquely determined from the coloring conditions at the crossings $c_{1}, c_{2}, \dots, c_{n-1}$ via the \textit{strand over-pass relations} \eqref{eq:over-first}--\eqref{eq:over-last}.
    When $(x_{1}, \dots, x_{n}) \in \mathcal{P}_{(n)}$, the coloring condition at $c_{n}$ is satisfied due to the \textit{strand over-pass relation} \eqref{eq:over-last2}:
    \[
    [x_{n+1},x_{n},y_{n}] = [x_{1}, x_{n}, y_{n}] = L(x_{1},x_{n})\cdots L(x_{2},x_{1})(y_{1}) = y_{1}.
    \]
    Assuming $y_{n+1} = y_{1}$, the $n$-tuple $(y_{1}, y_{2}, \dots, y_{n})$ satisfies the property $\mathcal{L}_{(n)}$ by Proposition~\ref{prop-RL-identities}~(5):
    \begin{align*}
        & L(y_{1},y_{n})L(y_{n},y_{n-1})\cdots L(y_{2},y_{1}) \\
        &= (L(y_{1},x_{1})L(x_{1},x_{n})L(x_{n},y_{n}))(L(y_{n},x_{n})L(x_{n},x_{n-1})L(x_{n-1},y_{n-1}))\cdots \\
        & \quad \cdots (L(y_{3},x_{3})L(x_{3},x_{2})L(x_{2},y_{2}))(L(y_{2},x_{2})L(x_{2},x_{1})L(x_{1},y_{1})) \\
        &= L(y_{1},x_{1})(L(x_{1},x_{n})L(x_{n},x_{n-1})\cdots L(x_{2},x_{1}))L(x_{1},y_{1})\\
        &= L(y_{1},x_{1})\mathrm{id}_{X}L(x_{1},y_{1}) \\
        &= \mathrm{id}_{X}.
    \end{align*}
    Similarly, $(y_{1}, y_{2}, \dots, y_{n}) \in\mathcal{R}_{(n)}$ by Proposition~\ref{prop-RL-identities}~(6).
    Thus, we find $(y_{1}, y_{2}, \dots, y_{n}) \in \mathcal{P}_{(n)}$.

    
    Considering the local coloring correspondence $g$ depicted in Figure~\ref{fig:strand_over_vertex}, the sequence \\
    $y_{1}, y_{2}, \dots, y_{n}$ satisfies the \textit{strand over-pass relations} \eqref{eq:over-first}--\eqref{eq:over-last} and $y_{n+1} = y_{1}$.
    Suppose $(y_{1},y_{2},\dots, y_{n}) \in \mathcal{P}_{(n)}$. From the \textit{strand over-pass relation} \eqref{eq:over-xn1}, 
    \[
        x_{n+1} = [y_{1},y_{n},x_{n}] = L(y_{1},y_{n})\cdots L(y_{2},y_{1})(x_{1}) = \mathrm{id}_{X}(x_{1}) = x_{1}.
    \]

    In addition, the sequence $x_{1},x_{2}, \dots, x_{n}$ satisfies the \textit{strand over-pass relations} \eqref{eq:over-x2}--\eqref{eq:over-xn1}, 
    and $(x_{1},x_{2},\dots,x_{n})\in \mathcal{R}_{(n)}$ by Proposition~\ref{prop-RL-identities}~(6):
    \begin{align*}
        & R(x_{1},x_{2})R(x_{2},x_{3})\cdots R(x_{n-1},x_{n})R(x_{n},x_{1}) \\
        &= (R(x_{1},y_{1})R(y_1,y_2)R(y_{2},x_{2}))(R(x_{2},y_{2})R(y_{2},y_{3})R(y_{3},x_{3}))\cdots \\
        & \quad \cdots (R(x_{n-1},y_{n-1})R(y_{n-1},y_{n})R(y_{n},x_{n}))(R(x_{n},y_{n})R(y_{n},y_{1})R(y_{1},x_{1})) \\
        &= R(x_{1},y_{1})(R(y_{1},y_{2})R(y_{2},y_{3})\cdots R(y_{n-1}, y_{n})R(y_{n},y_{1}))R(y_{1},x_{1}) \\
        &= R(x_{1},y_{1})\mathrm{id}_{X}R(y_{1},x_{1}) \\
        &= \mathrm{id}_{X}.
    \end{align*}
    Similarly, $(x_{1},x_{2}, \dots, x_{n}) \in \mathcal{L}_{(n)}$ by Proposition~\ref{prop-RL-identities}~(5).
    Thus, we find $(x_{1},x_{2}, \dots, x_{n}) \in \mathcal{P}_{(n)}$, and $g$ is well-defined.
    The coloring correspondences $f$ and $g$ are inverses of each other under the RIV move depicted in Figure~\ref{fig:strand_over_vertex}.

\begin{figure}[htbp]
        \centering
        \includegraphics{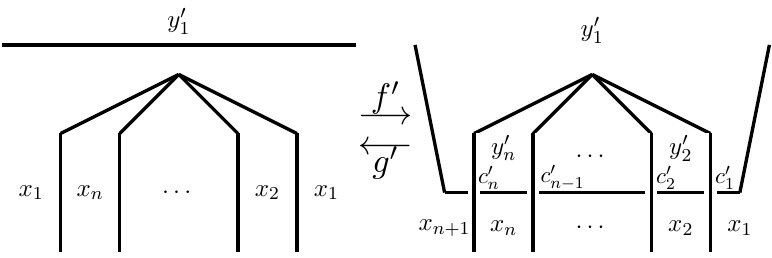}
        \caption{Colorings under the spatial graph move RIV (the strand passes under the vertex).}
        \label{fig:strand_under_vertex}
    \end{figure}

    For the other RIV move, depicted in Figure~\ref{fig:strand_under_vertex}, 
    the maps $f'$ 
    and $g'$ are similarly well-defined by Proposition~\ref{prop-RL-identities}~(7) and (8) together with the \textit{strand under-pass relations} \eqref{eq:under-first}--\eqref{eq:under-last}.
    In particular, for the map $g'$, if $(y_{1}', \dots, y_{n}')$ satisfies $\mathcal{R}_{(n)}$, then from the \textit{strand under-pass relation} \eqref{eq:under-last}, we have
    \[
        x_{n+1} = (x_{1})R(y_{1}',y_{2}')\cdots R(y_{n}',y_{1}') = x_{1}.
    \]
    Consequently, this coloring correspondence is bijective under the RIV move.
    
    \textbf{RV moves}

    The coloring correspondences under the RV moves are verified by Proposition~\ref{prop-RL-identities}~(1)--(4).

    \begin{figure}[htbp]
        \centering
        \includegraphics[width=\linewidth]{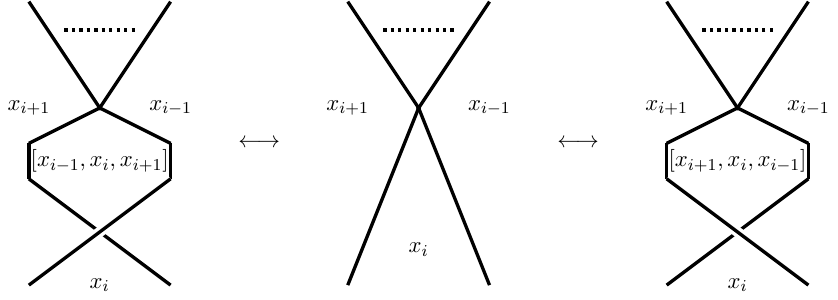}
        \caption{Colorings under the spatial graph move RV.}
        \label{fig:rv_invariance}
    \end{figure}

    Under the moves, a new region appears or vanishes. 
    The color of this area is determined uniquely by other values 
    as shown in Figure~\ref{fig:rv_invariance}.
    By Proposition~\ref{prop-RL-identities}~(1)--(4),
    \begin{align*}
        R(x_{i-1},[x_{i-1},x_{i},x_{i+1}])R([x_{i-1},x_{i},x_{i+1}],x_{i+1}) &= R(x_{i-1},x_{i})R(x_{i},x_{i+1}), \\
        L(x_{i+1},[x_{i-1},x_{i},x_{i+1}])L([x_{i-1},x_{i},x_{i+1}],x_{i-1}) &= L(x_{i+1},x_{i})L(x_{i},x_{i-1}), \\
        R(x_{i-1},[x_{i+1},x_{i},x_{i-1}])R([x_{i+1},x_{i},x_{i-1}],x_{i+1}) &= R(x_{i-1},x_{i})R(x_{i},x_{i+1}), \\
        L(x_{i+1},[x_{i+1},x_{i},x_{i-1}])L([x_{i+1},x_{i},x_{i-1}],x_{i-1}) &= L(x_{i+1},x_{i})L(x_{i},x_{i-1}). 
    \end{align*}

        These equations imply that the modified sequences of colors \\
        $(x_1,\dots,x_{i-1}, [x_{i-1},x_i,x_{i+1}],x_{i+1},\dots,x_n)$ and 
        $(x_1,\dots,x_{i-1}, [x_{i+1},x_i,x_{i-1}],x_{i+1},\dots,x_n)$\\ 
        satisfy the coloring condition around a vertex if and only if the original sequence \\
        $(x_1,\dots,x_{i-1}, x_i ,x_{i+1},\dots,x_n)$ does.

        Consequently, the coloring correspondences described above are well-defined and bijective.
\end{proof}


By Theorem~\ref{mainresult-1}, if $|X|<\infty$, the number of coloring maps is an invariant of spatial graphs. 

The concept of a KTQ palette was introduced in \cite{MR4259504}.
Here, we review the definition of an unoriented palette.
Let $(X,[\cdot,\cdot,\cdot]_{0})$ be an unoriented KTQ in the sense of \cite{MR4259504}.
For any $a,b,c\in X$, define the operators $H_{(a,b;2)}$ and $H_{(a,b;3)}$ by 
\[
    H_{(a,b;2)}(c) := [a,c,b]_{0} \quad \text{and} \quad H_{(a,b;3)}(c) := [a,b,c]_{0}.
\]

\begin{definition}[\cite{MR4259504}]\label{def-palette}
    An \textit{unoriented palette} $P$ of $(X,[\cdot,\cdot,\cdot]_{0})$ is a subset of $\bigcup_{n\geq 3}X^n$
    satisfying the following conditions:
    \begin{enumerate}[label=\upshape(\arabic*)]
        \item If $(a_{1}, a_{2}, \dots, a_{n}) \in P$, then $(a_{2}, \dots, a_{n}, a_{1}) \in P$.
        \item If $(a_{1}, a_{2}, \dots, a_{n}) \in P$, then it holds that
        \begin{align*}
            H_{(a_{n},a_{1};2)} \circ H_{(a_{n-1},a_{n};2)} \circ \cdots \circ H_{(a_{1},a_{2};2)} &= \mathrm{id}_{X},\\
            H_{(a_{n},a_{1};3)} \circ H_{(a_{n-1},a_{n};3)} \circ \cdots \circ H_{(a_{1},a_{2};3)} &= \mathrm{id}_{X}.
        \end{align*}
        \item If $(a_{1}, a_{2}, \dots, a_{n}) \in P$, then for any $a \in X$,
        \begin{align*}
            (a, H_{(a_{1},a_{2};2)}(a), H_{(a_{2},a_{3};2)} \circ H_{(a_{1},a_{2};2)}(a), \dots, H_{(a_{n-1},a_{n};2)} \circ \cdots \circ H_{(a_{1},a_{2};2)}(a)) &\in P, \\
            (a, H_{(a_{1},a_{2};3)}(a), H_{(a_{2},a_{3};3)} \circ H_{(a_{1},a_{2};3)}(a), \dots, H_{(a_{n-1},a_{n};3)} \circ \cdots \circ H_{(a_{1},a_{2};3)}(a)) &\in P.
        \end{align*}
        \item If $(a_{1}, a_{2}, \dots, a_{n}) \in P$,
        \begin{align*}
            (a_{1}, [a_{2},a_{3},a_{1}]_{0}, a_{3}, \dots, a_{n}) &\in P, \\
            (a_{1}, [a_{2},a_{1},a_{3}]_{0}, a_{3}, \dots, a_{n}) &\in P.
        \end{align*}
    \end{enumerate}
\end{definition}

These conditions are derived from the generalized Reidemeister moves RIV and RV for spatial graphs.
In particular, the original palette introduced in \cite{MR4259504} was defined as a subset of $\bigcup_{n \geq 1} X^n$.

In this paper, since we consider spatial graphs whose vertices have degree at least three, we restrict our attention to the case $n \geq 3$ for simplicity.
Let $(X,[\cdot,\cdot,\cdot])$ be an unoriented KTQ in our convention, and put $[a,b,c]_{0}=[b,a,c]$.
Note that we have the following operator equations:

\begin{align*}
    H_{(a,b;2)}(c) &= [a,c,b]_{0} = [c,a,b] = (c)R(a,b), \\
    H_{(a,b;3)}(c) &= [a,b,c]_{0} = [b,a,c] = L(b,a)(c).
\end{align*}

Using the equations $H_{(a,b;2)}(c) = (c)R(a,b)$ and $H_{(a,b;3)}(c) = L(b,a)(c)$ for $c\in X$, 
we can rewrite the compositions. 
For example, in the case $n=3$, for any $z\in X$:
\begin{align*}
        H_{(a_{3},a_{1};2)}\circ H_{(a_{2},a_{3};2)}\circ H_{(a_{1},a_{2};2)}(z)
        &= [a_{3},[a_{2},[a_{1},z,a_{2}]_{0},a_{3}]_{0},a_{1}]_{0} \\
        &= [[[z,a_{1},a_{2}],a_{2},a_{3}],a_{3},a_{1}] \\
        &= (z)R(a_{1},a_{2})R(a_{2},a_{3})R(a_{3},a_{1}), \\[1em]
        H_{(a_{3},a_{1};3)}\circ H_{(a_{2},a_{3};3)}\circ H_{(a_{1},a_{2};3)}(z)
        &= [a_{3},a_{1},[a_{2},a_{3},[a_{1},a_{2},z]_{0}]_{0}]_{0} \\
        &= [a_{1},a_{3},[a_{3},a_{2},[a_{2},a_{1},z]]] \\
        &= L(a_{1},a_{3})L(a_{3},a_{2})L(a_{2},a_{1})(z).
\end{align*}
Thus, the second axiom of an unoriented palette is exactly equivalent to the conditions $\mathcal{R}_{(n)}$ and $\mathcal{L}_{(n)}$. 
In other words, any unoriented palette $P$ must be a subset of $\bigcup_{n\geq 3}\mathcal{P}_{(n)}$.
Here we define the set $\mathcal{P}(X) := \bigcup_{n\geq 3}\mathcal{P}_{(n)}$.

\begin{proposition}\label{prop-maximality}
    Let $(X,[\cdot,\cdot,\cdot])$ be an unoriented KTQ in our convention, 
    and let \\$(X,[\cdot,\cdot,\cdot]_{0})$ be the corresponding unoriented KTQ in the convention of \cite{MR4259504} defined by $[a,b,c]_{0}=[b,a,c]$.
    Then the set $\mathcal{P}(X)$ is an unoriented palette of $(X, [\cdot,\cdot,\cdot]_0)$ in the sense of Definition~\ref{def-palette}.
    Moreover, it is the largest such palette with respect to inclusion. 
\end{proposition}
\begin{proof}
    To prove that $\mathcal{P}(X)$ is an unoriented palette, we must verify the four conditions of Definition~\ref{def-palette}. 

    Condition (2) holds by the definition of $\mathcal{P}_{(n)}$, and condition (1) follows immediately from Remark~\ref{remark-cyclic-ordering}. 
    Furthermore, conditions (3) and (4) are satisfied from the invariance under the RIV and RV moves, respectively, 
    as demonstrated in the proof of Theorem~\ref{mainresult-1}. 
    Thus, $\mathcal{P}(X)$ forms an unoriented palette. 

    Since every unoriented palette is contained in $\mathcal{P}(X)$, 
    $\mathcal{P}(X)$ is the largest unoriented palette.
\end{proof}


We can verify that some unoriented KTQs satisfy $\mathcal{P}_{(n)}=X^n$.

\begin{proposition}\label{prop-dehn-max}
    Let $(X, [\cdot, \cdot, \cdot])$ be a KTQ with the \textit{Dehn tribracket}.
    Then, $\mathcal{P}_{(n)}=X^n$ for any $n\geq 3$.
\end{proposition}
\begin{proof}
    It suffices to show that any $n$-tuple $(x_1, x_2, \dots, x_n) \in X^n$ satisfies the conditions $\mathcal{R}_{(n)}$ and $\mathcal{L}_{(n)}$.
    For any $z, a, b, c \in X$, we have
    \begin{align*}
        (z)R(a,b)R(b,c) &= (z a^{-1} b) b^{-1} c = z a^{-1} c = (z)R(a,c), \\
        L(c,b)L(b,a)(z) &= c b^{-1} (b a^{-1} z) = c a^{-1} z = L(c,a)(z).
    \end{align*}
    Thus, we obtain
    \begin{align*}
        (z)R(x_1,x_2)R(x_2,x_3)\cdots R(x_n,x_1) 
        &= (z)R(x_1,x_1) = z x_1^{-1} x_1 = z, \\
        L(x_1,x_n)L(x_n,x_{n-1})\cdots L(x_2,x_1)(z)
        &= L(x_1,x_1)(z) = x_1 x_1^{-1} z = z.
    \end{align*}
\end{proof}

\begin{proposition}\label{prop-unorialex-max}
Let $(X, [\cdot,\cdot,\cdot])$ be a KTQ with the unoriented Alexander tribracket (where the parameters $t$ and $s$ are involutive, i.e., $t^2 = s^2 = 1$). 
Then, $\mathcal{P}_{(2n)} = X^{2n}$ for any $n \geq 2$. 
Moreover, for any odd integer $2n+1 \geq 3$, 
we have $\mathcal{P}_{(2n+1)} = X^{2n+1}$ if $t = s = 1$, 
and $\mathcal{P}_{(2n+1)} = \emptyset$ otherwise.
\end{proposition}

\begin{proof}
For any $z, a, b, c \in X$, we evaluate the composition of the right operators:
\begin{align*}
(z)R(a, b)R(b, c) &= t(tz - tsa + sb) - tsb + sc \\
&= z - sa + sc.
\end{align*}
Similarly, for the left operators, we have:
\begin{align*}
L(c, b)L(b, a)(z) &= tc - tsb + s(tb - tsa + sz) \\
&= tc - ta + z.
\end{align*}
Thus, applying these compositions pairwise to the cyclic products of length $2n$, we get:
\begin{align*}
&(z)R(x_1, x_2)R(x_2, x_3)\cdots R(x_{2n-1}, x_{2n})R(x_{2n}, x_1) \\
&= z + \sum_{i=1}^{n-1}(-sx_{2i-1} + sx_{2i+1}) + (-sx_{2n-1} + sx_1) \\
&= z,
\end{align*}
and for the left operators:
\begin{align*}
&L(x_1, x_{2n})L(x_{2n}, x_{2n-1})\cdots L(x_3, x_2)L(x_2, x_1)(z) \\
&= (tx_1 - tx_{2n-1}) + \sum_{i=1}^{n-1}(tx_{2i+1} - tx_{2i-1}) + z \\
&= z.
\end{align*}
Hence, $\mathcal{P}_{(2n)} = X^{2n}$.

For an odd integer $2n+1 \geq 3$ and any $(x_1, \dots, x_{2n+1}) \in X^{2n+1}$, for any $z \in X$, we have
\begin{align*}
    &(z)R(x_1, x_2)R(x_2, x_3)\cdots R(x_{2n}, x_{2n+1})R(x_{2n+1}, x_1) \\ 
    &= \left(z + \sum_{i=1}^{n}(-sx_{2i-1} + sx_{2i+1})\right)R(x_{2n+1}, x_1) \\
    &= t(z - sx_1 + sx_{2n+1}) - tsx_{2n+1} + sx_1 \\
    &= tz + s(1 - t)x_1, \\  
    &L(x_1, x_{2n+1})L(x_{2n+1}, x_{2n})\cdots L(x_3, x_2)L(x_2, x_1)(z) \\
    &= L(x_1, x_{2n+1})\left(\sum_{i=1}^{n}(tx_{2i+1} - tx_{2i-1}) + z \right) \\
    &= tx_1 - tsx_{2n+1} + s(tx_{2n+1} - tx_1 + z) \\
    &= t(1 - s)x_1 + sz.
\end{align*}
If $\mathcal{P}_{(2n+1)} \neq \emptyset$, there exists a tuple for which both formulas equal $z$ for all $z \in X$. 
Taking the difference between the values for $z = 1$ and $z = 0$ yields $t = s = 1$, 
which implies $[x, y, z] = x - y + z$ and hence $\mathcal{P}_{(2n+1)} = X^{2n+1}$.
Otherwise, $\mathcal{P}_{(2n+1)}=\emptyset$.
\end{proof}



With the aim of finding examples distinct from the Dehn and unoriented Alexander tribrackets, 
we conducted a computer search in Python using the Z3 SMT solver. 
The unoriented KTQ axioms and the conditions 
\[ 
    \varnothing \neq \mathcal{P}_{(3)} \subsetneq X^3 
\] 
were encoded as constraints, and the solver was run with randomly selected seeds. 
As a result, we found an unoriented KTQ that is neither a Dehn tribracket nor an unoriented Alexander tribracket. 
The following example presents its ternary operation table and the corresponding set $\mathcal{P}_{(3)}$.

\begin{example}\label{ex-tribracket}
    Let $X$ be a finite set $\{1,2,3,4,5,6\}$, and define a tribracket operation by the following six matrices:
    \[
        \left[
        \begin{pmatrix}
        1 & 2 & 3 & 4 & 5 & 6 \\
        2 & 1 & 6 & 3 & 4 & 5 \\
        3 & 4 & 1 & 5 & 6 & 2 \\
        6 & 3 & 5 & 1 & 2 & 4 \\
        5 & 6 & 4 & 2 & 1 & 3 \\
        4 & 5 & 2 & 6 & 3 & 1
        \end{pmatrix},
        \begin{pmatrix}
        2 & 1 & 4 & 5 & 6 & 3 \\
        1 & 2 & 3 & 4 & 5 & 6 \\
        4 & 5 & 2 & 6 & 3 & 1 \\
        3 & 4 & 6 & 2 & 1 & 5 \\
        6 & 3 & 5 & 1 & 2 & 4 \\
        5 & 6 & 1 & 3 & 4 & 2
        \end{pmatrix},
        \begin{pmatrix}
        3 & 6 & 1 & 2 & 4 & 5 \\
        6 & 3 & 5 & 1 & 2 & 4 \\
        1 & 2 & 3 & 4 & 5 & 6 \\
        5 & 1 & 4 & 3 & 6 & 2 \\
        4 & 5 & 2 & 6 & 3 & 1 \\
        2 & 4 & 6 & 5 & 1 & 3
        \end{pmatrix},
        \right.
        \]

        \[
        \left.
        \begin{pmatrix}
        4 & 5 & 2 & 6 & 3 & 1 \\
        5 & 4 & 1 & 2 & 6 & 3 \\
        2 & 6 & 4 & 3 & 1 & 5 \\
        1 & 2 & 3 & 4 & 5 & 6 \\
        3 & 1 & 6 & 5 & 4 & 2 \\
        6 & 3 & 5 & 1 & 2 & 4
        \end{pmatrix},
        \begin{pmatrix}
        5 & 4 & 6 & 3 & 1 & 2 \\
        4 & 5 & 2 & 6 & 3 & 1 \\
        6 & 3 & 5 & 1 & 2 & 4 \\
        2 & 6 & 1 & 5 & 4 & 3 \\
        1 & 2 & 3 & 4 & 5 & 6 \\
        3 & 1 & 4 & 2 & 6 & 5
        \end{pmatrix},
        \begin{pmatrix}
        6 & 3 & 5 & 1 & 2 & 4 \\
        3 & 6 & 4 & 5 & 1 & 2 \\
        5 & 1 & 6 & 2 & 4 & 3 \\
        4 & 5 & 2 & 6 & 3 & 1 \\
        2 & 4 & 1 & 3 & 6 & 5 \\
        1 & 2 & 3 & 4 & 5 & 6
        \end{pmatrix}
        \right].
    \]
    
    Here, the value of the tribracket $[a,b,c]$ ($a,b,c \in X$) is given by the entry in the $b$-th row and the $c$-th column of the $a$-th matrix. 
    For example, $[2,6,4] = 3$ and $[4,3,1] = 2$.
    Direct computation gives the following $108$ elements of $\mathcal{P}_{(3)}$:
    
    \[
    \begin{aligned}
        \mathcal{P}_{(3)}=\left\{ \right.
        & (1, 1, 1), (1, 1, 2), (1, 1, 3), (1, 1, 4), (1, 1, 5), (1, 1, 6), (1, 2, 1), (1, 2, 2), (1, 3, 1), \\
        & (1, 3, 3), (1, 4, 1), (1, 4, 4), (1, 4, 6), (1, 5, 1), (1, 5, 5), (1, 6, 1), (1, 6, 4), (1, 6, 6), \\
        & (2, 1, 1), (2, 1, 2), (2, 2, 1), (2, 2, 2), (2, 2, 3), (2, 2, 4), (2, 2, 5), (2, 2, 6), (2, 3, 2), \\
        & (2, 3, 3), (2, 3, 5), (2, 4, 2), (2, 4, 4), (2, 5, 2), (2, 5, 3), (2, 5, 5), (2, 6, 2), (2, 6, 6), \\
        & (3, 1, 1), (3, 1, 3), (3, 2, 2), (3, 2, 3), (3, 2, 5), (3, 3, 1), (3, 3, 2), (3, 3, 3), (3, 3, 4), \\
        & (3, 3, 5), (3, 3, 6), (3, 4, 3), (3, 4, 4), (3, 5, 2), (3, 5, 3), (3, 5, 5), (3, 6, 3), (3, 6, 6), \\
        & (4, 1, 1), (4, 1, 4), (4, 1, 6), (4, 2, 2), (4, 2, 4), (4, 3, 3), (4, 3, 4), (4, 4, 1), (4, 4, 2), \\
        & (4, 4, 3), (4, 4, 4), (4, 4, 5), (4, 4, 6), (4, 5, 4), (4, 5, 5), (4, 6, 1), (4, 6, 4), (4, 6, 6), \\
        & (5, 1, 1), (5, 1, 5), (5, 2, 2), (5, 2, 3), (5, 2, 5), (5, 3, 2), (5, 3, 3), (5, 3, 5), (5, 4, 4), \\
        & (5, 4, 5), (5, 5, 1), (5, 5, 2), (5, 5, 3), (5, 5, 4), (5, 5, 5), (5, 5, 6), (5, 6, 5), (5, 6, 6), \\
        & (6, 1, 1), (6, 1, 4), (6, 1, 6), (6, 2, 2), (6, 2, 6), (6, 3, 3), (6, 3, 6), (6, 4, 1), (6, 4, 4), \\
        & (6, 4, 6), (6, 5, 5), (6, 5, 6), (6, 6, 1), (6, 6, 2), (6, 6, 3), (6, 6, 4), (6, 6, 5), (6, 6, 6)
        \left.\right\}.
    \end{aligned}
    \]
\end{example}

In Example~\ref{ex-tribracket}, $\mathcal{P}_{(3)}$ is a non-empty proper subset of $X^3$ since $(1, 1, 1) \in \mathcal{P}_{(3)}$ and $(1, 2, 3) \notin \mathcal{P}_{(3)}$. 
Consequently, the tribracket presented in Example~\ref{ex-tribracket} is neither a Dehn tribracket nor an unoriented Alexander tribracket by Propositions~\ref{prop-dehn-max} and~\ref{prop-unorialex-max}.

\section{coloring restrictions around vertices}

As shown in the previous section, if $\mathcal{P}_{(n)}=X^n$,
the coloring condition at vertices imposes no restriction on an $n$-tuple of region colors.
We therefore investigate an additional coloring restriction at vertices that yields further coloring invariants in the cases considered below.

Oshiro and Oyamaguchi~\cite{MR4259504} introduced the concept of \textit{palettes}, which are local coloring restrictions around a vertex preserved under RIV and RV moves.
The palettes depend on the chosen KTQ and the degree of the vertex.
In particular, they completely classified the palettes for Dehn $p$-colorings, 
which correspond to a special case of KTQs equipped with the unoriented Alexander tribracket.

Motivated by their work, we introduce an additional coloring restriction at vertices for the Dehn tribracket and the unoriented Alexander tribracket.
This restriction defines a subset
\[
    \mathcal{Q}_{(n)}(X,[\cdot,\cdot,\cdot])\subset X^n,
\]
whose elements represent admissible colorings around a vertex.
We then show that the resulting restricted colorings are preserved under generalized Reidemeister moves and therefore yield further invariants of spatial graphs.

\begin{definition}\label{def-condi-Q}
Let $(X,[\cdot,\cdot,\cdot])$ be an unoriented KTQ and let
$n\geq 3$ be an integer.

Define
\[
\mathcal{Q}_{(n)}(X,[\cdot,\cdot,\cdot])
=
\left\{
(x_1,x_2,\ldots,x_n)\in X^n
\;\middle|\;
x_i=[x_{i-1},x_i,x_{i+1}]
\text{ for all } i \pmod n
\right\}.
\]

When the KTQ $(X,[\cdot,\cdot,\cdot])$ is clear from the context,
we simply write $\mathcal{Q}_{(n)}$.
\end{definition}

We define
\[
\mathrm{Col}_{X}^{\mathcal{Q}}(D):=
\left\{
    c\in\mathrm{Col}_{X}(D)
    \;\middle|\;
    (x_1,\ldots,x_n)\in\mathcal{Q}_{(n)}
    \text{ at every vertex}
\right\}.
\]

\begin{lemma}\label{lem-commutative}
Let $(X,[\cdot,\cdot,\cdot])$ be an unoriented KTQ.
For any $a,b,c\in X$, the following conditions are equivalent:
\[
    a=[b,a,c]
    \quad\text{and}\quad
    a=[c,a,b].
\]
\end{lemma}

\begin{proof}
Suppose that $a=[b,a,c]$.
Since $(X,[\cdot,\cdot,\cdot])$ is unoriented, we have
\[
    [c,a,b]=[c,[b,a,c],b]=a.
\]
The converse follows similarly.
\end{proof}

By Lemma~\ref{lem-commutative}, the condition defining
$\mathcal{Q}_{(n)}$ is independent of the direction of the cyclic ordering, up to a cyclic permutation.


\begin{lemma}\label{lem-Q-RV-invariance}
    Let $(X, [\cdot, \cdot, \cdot])$ be an unoriented KTQ. If $D_1$ and $D_2$ are spatial graph diagrams related by an RV move, 
    then there is a bijection between $\operatorname{Col}_{X}^{\mathcal{Q}}(D_{1})$ and $\operatorname{Col}_{X}^{\mathcal{Q}}(D_{2})$.
\end{lemma}
\begin{proof}
    Under the moves depicted in Figure~\ref{fig:rv_invariance}, there is a coloring correspondence between the $n$-tuples
    \[
    (x_1, \dots, x_{i-1}, x_{i}, x_{i+1}, \dots, x_{n})
    \]
    and 
    \[
    (x_1, \dots, x_{i-1}, [x_{i-1}, x_{i}, x_{i+1}], x_{i+1}, \dots, x_{n})
    \]
    or
    \[
    (x_1, \dots, x_{i-1}, [x_{i+1}, x_{i}, x_{i-1}], x_{i+1}, \dots, x_{n}).
    \]
    
    Suppose $[x_{i-1}, x_{i}, x_{i+1}] = x_{i}$. Then,
    \[
    (x_1, \dots, x_{i-1}, [x_{i-1}, x_{i}, x_{i+1}], x_{i+1}, \dots, x_{n}) = (x_1, \dots, x_{i-1}, x_{i}, x_{i+1}, \dots, x_{n}).
    \]
    
    Conversely, suppose $[x_{i-1}, [x_{i-1}, x_{i}, x_{i+1}], x_{i+1}] = [x_{i-1}, x_{i}, x_{i+1}]$. By applying the unoriented property, we obtain
    \[
    [x_{i+1}, [x_{i-1}, [x_{i-1}, x_{i}, x_{i+1}], x_{i+1}], x_{i-1}] = [x_{i-1}, x_{i}, x_{i+1}],
    \]
    and simultaneously,
    \[
    [x_{i+1}, [x_{i-1}, x_{i}, x_{i+1}], x_{i-1}] = x_i.
    \]
    Thus, $[x_{i-1}, x_{i}, x_{i+1}] = x_i$. 
    
    We conclude that if the $n$-tuple 
    \[
    (x_1, \dots, x_{i-1}, [x_{i-1}, x_{i}, x_{i+1}], x_{i+1}, \dots, x_{n}) \in \mathcal{Q}_{(n)},
    \]
    then the $n$-tuple 
    \[
    (x_1, \dots, x_{i-1}, x_{i}, x_{i+1}, \dots, x_{n}) \in \mathcal{Q}_{(n)}.
    \]
    
    Similarly, we can verify the correspondence between
    \[
    (x_{1}, \dots, x_{i-1}, [x_{i+1}, x_{i}, x_{i-1}], x_{i+1}, \dots, x_{n})
    \]
    and 
    \[
    (x_{1}, \dots, x_{i-1}, x_{i}, x_{i+1}, \dots, x_{n})
    \]
    by Lemma~\ref{lem-commutative}.
\end{proof}


By Lemma~\ref{lem-Q-RV-invariance}, when coloring with an unoriented KTQ, there is a bijection between the restricted colorings under RV moves.

Furthermore, for certain KTQs, the sets $\mathrm{Col}_{X}^{\mathcal Q}(D)$ are invariant under all generalized Reidemeister moves.

\begin{proposition}\label{prop-dehn-Q-invariance}
    Let $D_1$ and $D_2$ be diagrams of a spatial graph, 
    and let $(X,[\cdot , \cdot , \cdot])$ be a KTQ with the Dehn tribracket. 
    Then there is a bijection between $\operatorname{Col}_{X}^{\mathcal{Q}}(D_1)$ and $\operatorname{Col}_{X}^{\mathcal{Q}}(D_2)$.
\end{proposition}
\begin{proof}
    By Lemma~\ref{lem-Q-RV-invariance}, it suffices to show that the condition defining $\mathcal{Q}_{(n)}$ is preserved under the RIV moves.
    Let the $n$-tuples $(x_1, x_2, \dots , x_{n})$ and $(y_1, y_{2}, \dots , y_{n})$ correspond to each other under the move depicted in Figure~\ref{fig:strand_over_vertex}.
    Here, $y_i$ (for $i=2, 3, \dots, n$) are represented as follows:
    \begin{align*}
        y_2 &= L(x_2,x_1)(y_{1}) = x_2 x_1^{-1}y_{1}, \\
        y_3 &= L(x_3,x_2)L(x_2, x_1)(y_1) = L(x_3, x_1)(y_{1}) = x_3 x_1^{-1} y_{1}, \\
        &\hspace{0.5em}\vdots \\
        y_{n} &= L(x_{n}, x_{n-1})\cdots L(x_2, x_1)(y_{1}) = L(x_{n}, x_{1})(y_{1}) = x_{n}x_{1}^{-1}y_{1}.
    \end{align*}
    For convenience, all subscripts are read modulo $n$; in particular, $x_{0}=x_{n}$, $y_{0}=y_{n}$, $x_{n+1}=x_{1}$, and $y_{n+1}=y_{1}$.
    We can calculate the tribracket as follows:
    \begin{align*}
        [y_{i-1}, y_{i}, y_{i+1}] 
        &= y_{i-1}y_{i}^{-1}y_{i+1} \\
        &= (x_{i-1}x_{1}^{-1}y_{1})(x_{i}x_{1}^{-1}y_{1})^{-1}(x_{i+1}x_{1}^{-1}y_{1}) \\
        &= [x_{i-1}, x_{i}, x_{i+1}]x_{1}^{-1}y_{1}.
    \end{align*}
    Suppose $[x_{i-1},x_{i},x_{i+1}] = x_i$. Then
    \[
    [y_{i-1}, y_{i}, y_{i+1}] = x_{i}x_{1}^{-1}y_{1} = y_i.
    \]
    Conversely, suppose $[y_{i-1}, y_{i}, y_{i+1}] = y_{i}$. Then,
    \begin{align*}
    [x_{i-1}, x_{i}, x_{i+1}] 
    &= [y_{i-1}, y_{i}, y_{i+1}]y_{1}^{-1}x_{1} \\
    &= y_{i}y_{1}^{-1}x_{1} \\
    &= (x_{i}x_{1}^{-1}y_{1})y_{1}^{-1}x_{1} \\
    &= x_i.
    \end{align*}
    For the other RIV move (depicted in Figure~\ref{fig:strand_under_vertex}, where the strand passes under the vertex), 
    the $n$-tuples $(x_1, x_2, \dots, x_{n})$ and $(y_{1}', y_{2}', \dots, y_{n}')$ are related as follows:
    \[
        y_{i}' = (y_{1}')R(x_{1},x_{2})\cdots R(x_{i-1},x_{i}) = y_{1}' x_{1}^{-1}x_{i}.
    \]
    Thus, 
    \[
        [y_{i-1}', y_{i}' , y_{i+1}'] = y_{1}'x_{1}^{-1}[x_{i-1}, x_{i}, x_{i+1}].
    \]
    Suppose $[x_{i-1},x_{i},x_{i+1}] = x_{i}$. Then 
    \[
        [y_{i-1}', y_{i}' , y_{i+1}'] = y_{1}'x_{1}^{-1}x_{i} = y_{i}'.
    \]
    Conversely, suppose $[y_{i-1}', y_{i}' , y_{i+1}']=y_{i}'$. Then
    \[
        [x_{i-1},x_{i},x_{i+1}] = x_{1}y_{1}'^{-1}y_{i}' = x_{i}.
    \]
    Therefore, the defining condition of $\mathcal{Q}_{(n)}$ is preserved under both cases of the RIV move.
\end{proof}

\begin{proposition}\label{prop-alex-Q-invariance}
    Let $D_1$ and $D_2$ be diagrams of a spatial graph.
    Let $(X,[\cdot,\cdot,\cdot])$ be a KTQ with the unoriented Alexander tribracket
    (where $t^2=s^2=1$).
    Then there is a bijection between
    \[
        \mathrm{Col}_X^{\mathcal Q}(D_1)
        \quad\text{and}\quad
        \mathrm{Col}_X^{\mathcal Q}(D_2).
    \]
\end{proposition}

\begin{proof}
    By Lemma~\ref{lem-Q-RV-invariance}, it suffices to show that the defining condition of $\mathcal Q_{(n)}$ is preserved under the RIV moves.

    Let $n$ be the degree of the considered vertex. Let the $n$-tuples
    \[
        (x_1,x_2,\ldots,x_n)\quad and \quad
        (y_1,y_2,\ldots,y_n)
    \]
    correspond to each other under the move depicted in Figure~\ref{fig:strand_over_vertex}.

    Here, for $i=2,\ldots,n$,
    \[
        y_i=L(x_i,x_{i-1})\cdots L(x_2,x_1)(y_1)=
        \begin{cases}
            tx_i-tx_1+y_1,
            & \text{if $i$ is odd},\\[0.3em]
            tx_i-tsx_1+sy_1,
            & \text{if $i$ is even}.
        \end{cases}
    \]

    For convenience, all subscripts are read modulo $n$; in particular,
    \[
        x_0=x_n,\quad
        y_0=y_n,\quad
        x_{n+1}=x_1,\quad
        y_{n+1}=y_1.
    \]

    We calculate
    \[
        [y_{i-1},y_i,y_{i+1}]=
        \begin{cases}
            t[x_{i-1},x_i,x_{i+1}]-tx_1+y_1,
            & \text{if $i$ is odd},\\[0.3em]
            t[x_{i-1},x_i,x_{i+1}]-tsx_1+sy_1,
            & \text{if $i$ is even}.
        \end{cases}
    \]

    If
    \[
        [x_{i-1},x_i,x_{i+1}]=x_i,
    \]
    then
    \[
        [y_{i-1},y_i,y_{i+1}]=
        \begin{cases}
            tx_i-tx_1+y_1,
            & \text{if $i$ is odd},\\[0.3em]
            tx_i-tsx_1+sy_1,
            & \text{if $i$ is even},
        \end{cases}
    \]
    and hence
    \[
        [y_{i-1},y_i,y_{i+1}]=y_i.
    \]

    Conversely, suppose that
    \[
        [y_{i-1},y_i,y_{i+1}]=y_i.
    \]

    If $i$ is odd, then
    \begin{align*}
        [x_{i-1},x_i,x_{i+1}]
        &=t[y_{i-1},y_i,y_{i+1}]+x_1-ty_1 \\
        &=ty_i+x_1-ty_1 \\
        &=t(tx_i-tx_1+y_1)+x_1-ty_1 \\
        &=x_i.
    \end{align*}

    If $i$ is even, then
    \begin{align*}
        [x_{i-1},x_i,x_{i+1}]
        &=t[y_{i-1},y_i,y_{i+1}]+sx_1-tsy_1 \\
        &=ty_i+sx_1-tsy_1 \\
        &=t(tx_i-tsx_1+sy_1)+sx_1-tsy_1 \\
        &=x_i.
    \end{align*}

    For the other RIV move
    (Figure~\ref{fig:strand_under_vertex}),
    the colors
    $y_i'$ ($i=2,\ldots,n$)
    are given by
    \[
        y_i'=(y_1')R(x_1,x_2)\cdots R(x_{i-1},x_i)=
        \begin{cases}
            sx_i-sx_1+y_1',
            & \text{if $i$ is odd},\\[0.3em]
            sx_i-tsx_1+ty_1',
            & \text{if $i$ is even}.
        \end{cases}
    \]

    We obtain
    \[
        [y_{i-1}',y_i',y_{i+1}']=
        \begin{cases}
            s[x_{i-1},x_i,x_{i+1}]
            -sx_1+y_1',
            & \text{if $i$ is odd},\\[0.3em]
            s[x_{i-1},x_i,x_{i+1}]
            -tsx_1+ty_1',
            & \text{if $i$ is even}.
        \end{cases}
    \]

    If
    \[
        [x_{i-1},x_i,x_{i+1}]=x_i,
    \]
    then
    \[
        [y_{i-1}',y_i',y_{i+1}']=
        \begin{cases}
            sx_i-sx_1+y_1',
            & \text{if $i$ is odd},\\[0.3em]
            sx_i-tsx_1+ty_1',
            & \text{if $i$ is even},
        \end{cases}
    \]
    and therefore
    \[
        [y_{i-1}',y_i',y_{i+1}']=
        y_i'.
    \]

    Conversely, suppose that
    \[
        [y_{i-1}',y_i',y_{i+1}']=y_i'.
    \]

    If $i$ is odd, then
    \begin{align*}
        [x_{i-1},x_i,x_{i+1}]
        &=s[y_{i-1}',y_i',y_{i+1}']+x_1-sy_1' \\
        &=sy_i'+x_1-sy_1' \\
        &=s(sx_i-sx_1+y_1')+x_1-sy_1' \\
        &=x_i.
    \end{align*}

    If $i$ is even, then
    \begin{align*}
        [x_{i-1},x_i,x_{i+1}]
        &=s[y_{i-1}',y_i',y_{i+1}']+tx_1-tsy_1' \\
        &=sy_i'+tx_1-tsy_1' \\
        &=s(sx_i-tsx_1+ty_1')+tx_1-tsy_1' \\
        &=x_i.
    \end{align*}

    Therefore, the defining condition of $\mathcal{Q}_{(n)}$ is preserved under both types of RIV moves.

    Together with Lemma~\ref{lem-Q-RV-invariance}, this yields a bijection between
    $\mathrm{Col}_X^{\mathcal Q}(D_1)$ and $\mathrm{Col}_X^{\mathcal Q}(D_2)$.
\end{proof}


The following example illustrates that the restricted coloring numbers associated with $\mathcal{Q}_{(n)}$ can distinguish spatial graphs that are not distinguished by unrestricted coloring numbers.

A \emph{constituent knot} of a spatial graph is the image of a cycle (a subgraph homeomorphic to $S^1$) in the spatial graph.
It is known that the multiset of constituent knots of a spatial graph is an invariant \cite{MR946218}.

\begin{example}\label{EX-application}
    The two spatial graphs depicted in Figure~\ref{fig:theta-curves} are embedded $\theta$-curves taken from the table of $\theta$-curves in \cite{MR2507922}.
    Each of them has one trefoil knot and two trivial knots as its constituent knots.
    
    \begin{figure}[htbp]
        \centering
        \begin{subfigure}{0.45\textwidth}
            \centering
            \includegraphics[width=\linewidth]{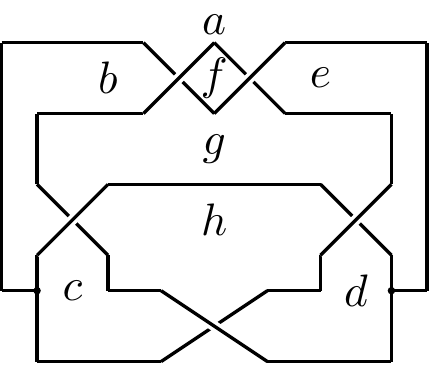}
            \caption{The $\theta$-curve $5_2$}
            \label{fig:theta52}
        \end{subfigure}
        \hfill
        \begin{subfigure}{0.45\textwidth}
            \centering
            \includegraphics[width=\linewidth]{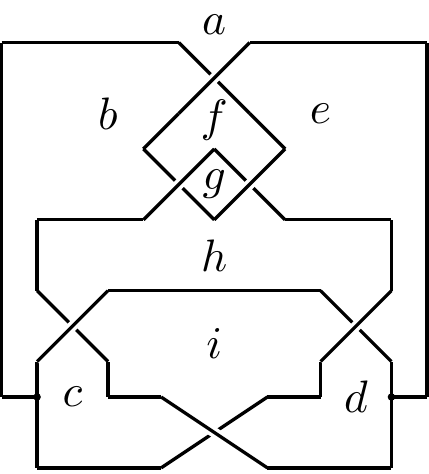}
            \caption{The $\theta$-curve $6_2$}
            \label{fig:theta62}
        \end{subfigure}
        \caption{Two $\theta$-curves used in this example.}
        \label{fig:theta-curves}
    \end{figure}

    Let $A_{4}$ be the alternating group of degree $4$, and let $(A_{4}, [\cdot,\cdot,\cdot])$ be the Dehn tribracket on $A_{4}$.
    The diagrams depicted in Figure~\ref{fig:theta-curves} have regions $\mathcal{R}(5_2)=\{a,b,c,d,e,f,g,h\}$ and $\mathcal{R}(6_2)=\{a,b,c,d,e,f,g,h,i\}$, respectively.
    We compute the coloring numbers of these two curves by the Dehn tribracket:
    \begin{align*}
        \mathrm{Col}_{(A_{4},[\cdot,\cdot,\cdot])}(5_2) 
        &= \bigl\{
            C:\mathcal{R}(5_2)\to A_{4} \mid 
            C(f)=[C(a),C(b),C(g)]=[C(g),C(e),C(a)], \\
        &\qquad C(h)=[C(g),C(b),C(c)]=[C(d),C(e),C(g)]=[C(c),C(a),C(d)]
        \bigr\}, \\
        \mathrm{Col}_{(A_{4},[\cdot,\cdot,\cdot])}(6_2) 
        &= \bigl\{
            C:\mathcal{R}(6_2)\to A_{4} \mid 
            C(f)=[C(b),C(a),C(e)], \\
        &\qquad C(g)=[C(f),C(b),C(h)]=[C(h),C(e),C(f)], \\
        &\qquad C(i)=[C(h),C(b),C(c)]=[C(d),C(e),C(h)]=[C(c),C(a),C(d)]
        \bigr\}. 
    \end{align*}
    For $D \in \{5_2, 6_2\}$,
    \begin{align*}
        \mathrm{Col}_{(A_{4},[\cdot,\cdot,\cdot])}^{\mathcal{Q}}(D)
        &= \bigl\{
            C\in \mathrm{Col}_{(A_{4},[\cdot,\cdot,\cdot])}(D) \mid
            C(a)=[C(c),C(a),C(b)],\\
        &\qquad C(b)=[C(a),C(b),C(c)],\, C(c)=[C(b),C(c),C(a)],\\
        &\qquad C(a)=[C(e),C(a),C(d)],\, C(d)=[C(a),C(d),C(e)],\\
        &\qquad C(e)=[C(d),C(e),C(a)]
        \bigr\}.
    \end{align*}

    Using GAP~\cite{GAP}, we obtain the numbers of colorings:
    \[
        |\mathrm{Col}_{(A_{4},[\cdot,\cdot,\cdot])}(5_2)| = |\mathrm{Col}_{(A_{4},[\cdot,\cdot,\cdot])}(6_2)| = 2880,
    \]
    whereas
    \[
        |\mathrm{Col}_{(A_{4},[\cdot,\cdot,\cdot])}^{\mathcal{Q}}(5_2)| = 396 \neq 108 = |\mathrm{Col}_{(A_{4},[\cdot,\cdot,\cdot])}^{\mathcal{Q}}(6_2)|.
    \]

    In conclusion, these $\theta$-curves cannot be distinguished by their constituent knots or by the unrestricted coloring numbers $|\mathrm{Col}_{(A_{4},[\cdot,\cdot,\cdot])}(D)|$, whereas they are distinguished by the restricted coloring numbers $|\mathrm{Col}_{(A_{4},[\cdot,\cdot,\cdot])}^{\mathcal{Q}}(D)|$.
\end{example}


\section*{tool use}
The author used Python in Example~\ref{ex-tribracket} to search for a suitable operation table and to compute $\mathcal{P}_{(3)}$.

The author used GAP in Example~\ref{EX-application} to calculate the number of colorings.

The author also used Microsoft Copilot and Google Gemini to assist with English language editing, 
LaTeX formatting, the presentation of some mathematical arguments, 
and the writing and editing of the programming code used for the computer search. 
All mathematical statements, proofs, and codes were independently verified by the author, 
who takes full responsibility for the entire contents of the paper.

\section*{acknowledgments}
The author would like to express sincere gratitude to his supervisor Seiichi Kamada for his continuous guidance and invaluable advice. 

\bibliographystyle{alpha}
\bibliography{references}

\end{document}